\documentclass{amsart}
\usepackage{graphicx,amssymb,epsfig}
\newtheorem{theorem}{Theorem}[section]
\newtheorem{proposition}[theorem]{Proposition}
\newtheorem{lemma}[theorem]{Lemma}
\newtheorem{corollary}[theorem]{Corollary}
\newtheorem{remark}[theorem]{Remark}
\newtheorem{example}[theorem]{Example}

\newtheorem{definition}[theorem]{Definition}

\newcommand{\bth}{\begin{theorem}}
\newcommand{\bpr}{\begin{proposition}}
\newcommand{\epr}{\end{proposition}}
\newcommand{\bco}{\begin{corollary}}
\newcommand{\eco}{\end{corollary}}
\newcommand{\ble}{\begin{lemma}}
\newcommand{\ele}{\end{lemma}}
\newcommand{\bre}{\begin{remark}\rm}
\newcommand{\ere}{\end{remark}}
\newcommand{\bex}{\begin{example}\rm}
\newcommand{\eex}{\end{example}}
\newcommand{\bde}{\begin{definition}\rm}
\newcommand{\ede}{\end{definition}}

\def\lrar{{\ra 2}}

\def\tensor{\otimes}
\def\sp#1{{SP}^{#1}}
\def\bsp#1{\overline{{SP}}^{#1}}
\def\tsp#1{\widetilde{{SP}}^{#1}}
\def\spy{{SP}^{\infty}}

\def\map#1{\hbox{Map}_{#1}}
\def\bmap#1{\hbox{Map}^*_{#1}}
\def\tp#1{\hbox{TP}^{#1}}

\usepackage[all]{xy}

\def\la#1{\hbox to #1pc{\leftarrowfill}}
\def\ra#1{\hbox to #1pc{\rightarrowfill}}
\def\fract#1#2{\raise4pt\hbox{$ #1 \atop #2 $}}
\def\decdnar#1{\phantom{\hbox{$\scriptstyle{#1}$}}
\left\downarrow\vbox{\vskip10pt\hbox{$\scriptstyle{#1}$}}\right.}

\def\lrar{{\ra 2}}

\def\tensor{\otimes}
\def\cohdim{\hbox{cohdim}}

\def\bbz{{\mathbb Z}}
\def\bbf{{\mathbb F}}
\def\bbp{{\mathbb P}}
\def\bbr{{\mathbb R}}

\def\bbc{{\mathbb C}}
\def\bbr{{\mathbb R}}

\def\bbn{{\mathbb N}}

\def\map#1{\hbox{Map}_{#1}}
\def\bmap#1{\hbox{Map}^*_{#1}}

\def\tensor{\otimes}

\begin{document}

\title{Symmetric Products, Duality and Homological Dimension of
Configuration Spaces}
\author{Sadok Kallel}
\address{Universit\'e des Sciences et Technologies de Lille\\
  Laboratoire Painlev\'e, U.F.R de Math\'ematiques\\
 59655 Villeneuve d'Ascq, France}
\email{sadok.kallel@math.univ-lille1.fr}

\maketitle
\centerline{To Fred Cohen on his 60th Birthday}

\begin{abstract} We discuss various aspects of "braid spaces'' or
  configuration spaces of unordered points on manifolds.  First we describe
  how the homology of these spaces is affected by puncturing the underlying
  manifold, hence extending some results of Fred Cohen, Goryunov and Napolitano.
  Next we obtain a precise bound for the cohomological dimension of braid
  spaces.  This is related to some sharp and useful connectivity bounds that
  we establish for the reduced symmetric products of any simplicial complex.
  Our methods are geometric and exploit a dual version of configuration spaces
  given in terms of truncated symmetric products.  We finally refine and then
  apply a theorem of McDuff on the homological connectivity of a map from
  braid spaces to some spaces of ``vector fields''.
\end{abstract}


\section{Introduction}

Braid spaces or configuration spaces of \textit{unordered pairwise
distinct} points on manifolds have important applications to a
number of areas of mathematics and physics. They were of crucial use
in the seventies in the work of Arnold on singularities and then
later in the eighties in work of Atiyah and Jones on instanton
spaces in gauge theory. In the nineties they entered in many works
on the homological stability of holomorphic mapping spaces. No more
important perhaps had been their use than in stable homotopy theory
in the sixties and early seventies through the work of Milgram, May,
Segal and Fred Cohen who worked out the precise connection with loop
space theory. This work has led in particular to the proof of
Nishida's nilpotence theorem and to Mahowald's infinite family in
the stable homotopy groups of spheres to name a few.

Given a space $M$, define $B(M,n)$ to be the space of finite subsets
of $M$ of cardinality $n$.  This is usually referred to as the
$n$-th "braid space" of $M$ and in the literature it is often
denoted by $C_n(M)$ \cite{atiyah,bct,cohen}. Its fundamental group
written $Br_n(M)$ is the ``braid group'' of $M$. The object of this
paper is to study the homology of braid spaces and the main approach
we adopt is that of duality with the symmetric products. In so doing
we take the opportunity to refine and elaborate on some classical
material. Next is a brief content summary.

Section \ref{braids} describes the homotopy type of braid spaces of
some familiar spaces and discusses orientation issues.  Section
\ref{tp} introduces truncated products as in \cite{bcm,lm}, states
the duality with braid spaces and then proves our first main result
on the cohomological dimension of braid spaces. Section
\ref{punctured} uses truncated product constructions to split in an
elementary fashion the homology of braid spaces for punctured
manifolds. In section \ref{bounds} we prove our sharp connectivity
result for \textit{reduced} symmetric products of CW complexes which
seems to be new and a significant improvement on work of Nakaoka and
Welcher \cite{welcher}. In section \ref{spec} we make the link
between the homology of symmetric and truncated products by
discussing a spectral sequence introduced by C.F. Bodigheimer, Fred
Cohen and R.J. Milgram and exploited by them to study ``braid
homology'' $H_*(B(M,n))$. Finally Section \ref{stability} completes
a left out piece from McDuff and Segal's work on configuration
spaces \cite{dusa}.  In that paper, $H_*(B(M,n))$, for closed
manifolds $M$, is compared to the homology of some spaces of
"compactly supported vector fields" on $M$ and the main theorem
there states that these homologies are isomorphic up to a range that
increases with $n$. We make this range more explicit and use it for
example to determine the abelianization of the braid groups of a
closed Riemann surface. A final appendix collects some homotopy
theoretic properties of section spaces that we use throughout.

Below are precise statements of our main results which we have
divided up into three main parts.  Unless explicitly stated, all
spaces are assumed to be connected. The $n$-th symmetric group is
written $\mathfrak{S}_n$.

\subsection{Connectivity and Cohomological Dimension}
For $M$ a manifold, write $H^*(M,\pm \bbz )$ the cohomology of $M$
with coefficients in the orientation sheaf $\pm\bbz$; that is
$H^*(M,\pm\bbz)$ is the homology of $Hom_{\bbz[\pi_1(X)]}(C_*(\tilde
M), \bbz)$, where $C_*(\tilde M)$ is the singular chain complex of
the universal cover $\tilde M$ of $M$, and where the action of (the
class of) a loop on the integers $\bbz$ is multiplication by $\pm 1$
according to whether this loop preserves or reverses orientation.
Similarly one defines $H_*(M,\pm\bbz ):= H_*(C_*(\tilde
M)\tensor_{\bbz[\pi_1(x)]}\bbz)$.

\bre\label{twisted} (see Lemma \ref{folklore}) When $M$ is simply
connected
  and $\dim M:=d > 2$, $\pi_1(B(M,k))=\mathfrak{S}_k$ and $\tilde B(M,k) =
  F(M,k)\subset M^k$ is the subspace of $k$ \textit{ordered} pairwise distinct
  points in $M$ (\S\ref{braids}). It follows that
  $H^*(B(M,k);\pm\bbz)$ is the homology of the chain complex
  $Hom_{\bbz[\mathfrak{S}_k]}(C_*(F(M,k),\bbz)$ where $\mathfrak{S}_k$ acts on
  $\bbz$ via $\sigma (1) = (-1)^{sg( \sigma)\cdot d }$ and $sg (\sigma
  )$ is the sign of the permutation $\sigma\in\mathfrak{S}_k$. \ere

We denote by $\cohdim_{\pm\bbz}(M)$ (cohomological dimension) the
smallest integer with the property that
$$H^i(M;\pm\bbz ) = 0\ \ ,\ \ \forall i > \cohdim_{\pm\bbz }(M)$$
If $M$ is orientable, then $H^*(M,\pm\bbz ) = H^*(M,\bbz )$ and
$\cohdim_{\pm\bbz}(M) = \cohdim (M)$; the cohomological dimension of
$M$.

A space $X$ is $r$-connected if $\pi_i(X)=0$ for $0\leq i\leq r$.
The connectivity of $X$; $conn(X)$, is the largest integer with such
a property.  This connectivity is infinite if $X$ is contractible.
The following is our first main result

\bth\label{main3} Let $M$ be a compact manifold of dimension $d\geq
1$,
  with boundary $\partial M$, and let $U\subset M$ be a closed subset such
  that $U\cap\partial M= \emptyset$ and $M-U$ connected.
We denote by $r$ the connectivity of $M$
  if $U\cup\partial M=\emptyset$, or the connectivity of the quotient
  $M/U\cup\partial M$ if $U\cup\partial M\neq \emptyset$. We assume $0\leq
  r<\infty$ and $k\geq 2$. Then
$$
\cohdim_{\pm\bbz}(B(M-U,k)) \leq
\begin{cases} (d-1)k-r+1, &\hbox{if}\ U\cup\partial M=\emptyset\\
(d-1)k-r, &\hbox{if}\ U\cup\partial M\neq \emptyset
\end{cases}
$$ When $M$ is even dimensional orientable, then replace
$\cohdim_{\pm\bbz}$ by $\cohdim$.
\end{theorem}

\bre\label{numbered} We check this theorem against some known
examples:
\begin{enumerate}
\item $B(S^d-\{p\},2)=B(\bbr^d,2)\simeq\bbr P^{d-1}$
(see Section \ref{tp}) and $\cohdim_{\pm\bbz}(B(\bbr^d,2)) =
2(d-1)-r = d-1 =\cohdim_{\pm\bbz}(\bbr P^{d-1})$ indeed, where
$r=d-1= conn(S^d)$.
\item
$B(S^d,2)\simeq\bbr P^d$ (see Section \ref{tp}) and
$\cohdim_{\pm\bbz} (B(S^d,2)) = d$\ in agreement with our formula.
\item
It is known that for odd primes $p$ and $d\geq 2$,
$H^{(d-1)(p-1)}(B(\bbr^d,p);\bbf_p)$ is non-trivial and an
isomorphic image of $H^{(d-1)(p-1)}(\mathfrak{S}_p;\bbf_p)$
\cite{ossa, vassiliev}. Our result states that, at least for even
$d$, no higher homology can occur. The cohomological dimension of
$B(\bbr^d,k)$ when using $\bbf_2$ coefficients is known to be
$(k-\alpha (k))\cdot (d-1)$ where $\alpha (k)$ is the number of 1's
in the dyadic decomposition of $k$ (see \cite{frido}). In the case
$d=2$, $B(\bbr^2,k)$ is the classifying space of Artin braid group
$B_k:= Br_k(\bbr^2)$ and is homotopy equivalent to a
$(k-1)$-dimensional CW complex so that $\cohdim (B(\bbr^2,k))\leq
k-1$ in agreement with our calculation.
\end{enumerate}
\ere

\bre The theorem applies to when $M=S^1$ and $U$ is either empty or
a single point. In that case $M-U=S^1,\bbr$. But one knows that for
$k\geq 1$, $B(S^1,k)\simeq S^1$ (Proposition \ref{cns1}) and
$B(\bbr,k)$ is contractible. \ere

\bco\label{twocomplexes} Let $S$ be a Riemann surface and
  $Q\subset S$ a finite subset.  Then $H^i(B(S-Q,k)) = 0$ if $i\geq k+1$ and
  $Q\cup\partial S\neq\emptyset$ ; or if $i> k+1$ and $Q\cup\partial
  S=\emptyset$.  \eco

This corollary gives an extension of the "finiteness'' result of
\cite{nap1}.  When $S$ is an open surface, then $B(S,k)$ is a Stein
variety and hence its homology vanishes above the complex dimension;
i.e. $H_i(B(S,k)) = 0$ for $i> k$. This also agrees with the above
computed bounds.

The proof of Theorem \ref{main3} relies on a useful connectivity
result of Nakaoka (Theorem \ref{nakak}). We also use this result to
produce sharp connectivity bounds for the \textit{reduced} symmetric
products \S\ref{bounds}.  Recall that $\sp{n}(X)$, the $n$-th
symmetric product of $X$, is the quotient of $X^n$ by the
permutation action of the symmetric group $\mathfrak{S}_n$ so that
$B(X,n)\subset\sp{n}(X)$ is the subset of configurations of distinct
points. We always assume $X$ is based so there is an embedding
$\sp{n-1}(X)\hookrightarrow\sp{n}(X)$ given by adjoining the
basepoint, with cofiber $\bsp{n}(X)$ the "$n$-th reduced symmetric"
product of $X$. The following result expresses the connectivity of
$\bsp{n}X$ in terms of the connectivity of $X$.

\bth\label{connectivity} Suppose $X$ is a based $r$-connected
simplicial
  complex with $r\geq 1$. Then $\bsp{n}(X)$ is $2n+r-2$- connected.
\end{theorem}

In particular the embedding $\sp{n-1}(X)\lrar\sp{n}(X)$ induces
homology isomorphisms in degrees up to $(2n+r-3)$.  The proof of
this theorem is totally inspired from \cite{kk} where similar
connectivity results are stated, and it uses the fact that the
homology of symmetric products only depends on the homology of the
underlying complex \cite{dold}.  Note that the bound $2n+r-2$ is
sharp as is illustrated by the case $X=S^2$, $r=1$ and
$\bsp{n}(S^2)=S^{2n}$.  A slightly weaker connectivity bound than
ours can be found in \cite{welcher}, corollary 4.9.

Note that Theorem \ref{connectivity} is stated for simply connected
spaces.  To get connectivity results for reduced symmetric products
of a compact Riemann surface for example we use geometric input from
\cite{ks2}. This applies to any two dimensional complex.

\bpr\label{conntwo} Let $X = \bigvee^wS^1\cup (D^2_1\cup\cdots\cup
  D^2_r)$ be a two dimensional CW complex with one skeleton a bouquet of $w$
  circles.  Then $\bsp{n}X$ is $(2n-min(w,n)-1)$-connected.
\epr


\subsection{Puncturing Manifolds}
We give generalizations and a proof simplification of results of
Napolitano \cite{nap1,nap2}.  For $S$ a two dimensional topological
surface, $p$ and the $p_i$ points in $S$, it was shown in
\cite{nap1} that for field coefficients $\bbf$,
\begin{equation}\label{secondsplit}
H^j(B(S -\{p_1,p_2\},n);\bbf )\cong \bigoplus_{t=0}^nH^{j-t}(B(S
-\{p\},n-t);\bbf )
\end{equation}
Here and throughout $H^* = 0$ when $*<0$ and $B(X,0)$ is basepoint.
When $S$ is a closed orientable surface and $\bbf=\bbf_2$,
\cite{nap1} establishes furthermore a splitting
\begin{equation}\label{firstsplit}
H^j(B(S,n);\bbf_2)\cong H^{j}(B(S -\{p\},n);\bbf_2)\oplus
H^{j-2}(B(S -\{p\},n-1);\bbf_2)
\end{equation}
Similar splittings occur in \cite{cohen2, goryunov}. These
splittings as we show extend to any closed topological manifold $M$
and to any number of punctures. If $V$ is a vector space, write
$V^{\oplus k}:= V\oplus\cdots\oplus V$ ($k$-times). Given positive
integers $r$ and $s$, we write $p(r,s)$ the number of ways we can
partition $s$ into a sum of $r$ \textit{ordered} positive (or null)
integers. For instance $p(1,s)=1$, $p(2,s)=s+1$ and $p(r,1)=r$.

\bth\label{main} Let $M$ be a closed connected manifold of
  dimension $d$ and $p\in M$. Then
\begin{equation}\label{main11}
H^j(B(M,n);\bbf_2)\cong H^{j}(B(M-\{p\},n);\bbf_2)\oplus
H^{j-d}(B(M-\{p\},n-1);\bbf_2)
\end{equation}
If moreover $M$ is oriented and even dimensional, then
\begin{equation}\label{main22}
H^j(B(M-\{p_1,\cdots, p_k\},n);\bbf )\cong \bigoplus_{0\leq r\leq n}
H^{j - (n-r)(d-1)} (B(M-\{p\},r);\bbf )^{\oplus p(k-1,n-r)}
\end{equation}
For an arbitrary closed manifold, (\ref{main22}) is still true with
$\bbf_2$-coefficients.
\end{theorem}

\bre As an example we can set $M=S^2, k=2=d$ and obtain the additive
splitting $H^j(B(\bbc^*,n);\bbf )\cong \bigoplus_{0\leq r\leq n}
H^{j - (n-r)} (B(\bbc ,r);\bbf )$\ as in (\ref{secondsplit}), where
$\bbc^*$ is the punctured disk (this isomorphism holds integrally
according to \cite{goryunov}).  Note that the left hand side is the
homology of the hyperplane arrangement of `Coxeter type'' $B_n$;
that is $B(\bbc^*,n)$ is an Eilenberg-MacLane space
$K(Br_n(\bbc^*),1)$ with fundamental group isomorphic to the
subgroup of Artin's braids $Br_{n+1}(\bbc )$ consisting of those
braids which leave the last strand fixed. It can be checked that the
abelianization of this group for $n\geq 2$ is $\bbz^2$ which is
consistent with the calculation of $H^1$ obtained from the above
splitting. \ere

Napolitano's approach to (\ref{secondsplit}) is through spectral
sequence arguments and ``resolution of singularities'' as in
Vassiliev theory.  Our approach relies on a simple geometric
manipulation of the truncated symmetric products as discussed
earlier (see Section \ref{punctured}). Theorem \ref{main} is a
consequence of combining a Poincar\'e-Lefshetz duality statement,
the identification of truncated products of the circle with real
projective space \cite{mostovoy} and a homological splitting result
due to Steenrod (Section \ref{tp}). Note that the splitting in (3)
is no longer true with coefficients other than $\bbf_2$ and is
replaced in general by a long exact sequence (lemma
\ref{longexact}).

\subsection{Homological Stability}
This is the third and last part of the paper.
 For $M$ a closed smooth manifold of dimension
$\dim M = d$, let $\tau^+M$ be the fiberwise one-point
compactification of the tangent bundle $\tau M$ of $M$ with fiber
$S^d$. We write $\Gamma (\tau^+M)$ the space of sections of
$\tau^+M$. Note that this space has a preferred section (given by
the points at infinity). There are now so called "scanning" maps for
any $k\in\bbn$ \cite{dusa, bct, quarterly}
\begin{equation}\label{firstscan}
S_k : B(M,k)\lrar \Gamma_k (\tau^+M )
\end{equation}
where $\Gamma_k(\tau^+M)$ is the component of degree $k$ sections
(see \S6.2). In important work, McDuff shows that $S_k$ induces a
homology isomorphism through a range that increases with $k$. In
many special cases, this range needs to be made explicit and this is
what we do next.

We say that a map $f: X\rightarrow Y$ is homologically $k$-connected
(or a homology equivalence up to degree $k$) if $f_*$ in homology is
an isomorphism up to and including degree $k$.

\bpr\label{main4} Let $M$ be a closed manifold of dimension $d\geq
2$. Assume the map $+:B(M-p,k)\lrar B(M-p,k+1)$ which consists of
adding a point near $p\in M$ (see Section \ref{stability}) is
homologically $s(k)$-connected. Then scanning $S_k$ is homologically
$s(k-1)$-connected. Moreover $s(k)\geq [k/2]$ (Arnold). \epr

When $k=1$, we give some information about $S_1 :
M\lrar\Gamma_1(\tau^+M)$ in lemma \ref{s1}. Note that $s(k)$ is an
increasing function of $k$. Arnold's inequality $s(k)\geq [k/2]$ is
proven in \cite{segal1}. This bound is far from being optimal in
some cases since for instance, for $M$ a compact Riemann surface,
$s(k) = k-1$ \cite{ks}. Note that the actual connectivity of the map
$+:B(M-p,k)\lrar B(M-p,k+1)$ is often $0$ since if $\dim M > 2$,
this map is never trivial on $\pi_1$ (see Lemma \ref{folklore}).

The utility of Proposition \ref{main4} is that in some particular
cases, knowledge of the homology of braid spaces in a certain range
informs on the homology of some mapping spaces. Here's an
interesting application to computing the abelianization of the braid
group of a surface (this was an open problem for some time).

\bco \label{cor1} For $S$ a compact Riemann surface of genus $g\geq
1$, and $k\geq 3$, we have the isomorphism: $H_1(B(S,k);\bbz ) =
\bbz_2\oplus\bbz^{2g}$. \eco

\begin{proof} $\tau^+S$ is trivial since $S$ is stably parallelizable
  and $\Gamma (\tau^+S) \simeq Map (S, S^2)$. Suppose $S$ has odd genus,
  then $S_k : H_1(B(S,k))\lrar H_1(Map_k(S,S^2))$ is degree preserving
(where degree is $k$)
  and according to Proposition \ref{main4} it is an isomorphism when
  $k\geq 3$ using the bound provided by Arnold.
  But $\pi:=\pi_1(Map_k(S,S^2))$ was
  computed in \cite{contemp} and it is some extension
  $$0\lrar\bbz_{2|k|}\lrar\pi\lrar \bbz^{2g}\lrar 0$$
  with a generator $\tau$ and torsion free generators $e_1,\ldots, e_{2g}$
  with non-zero commutators $[e_i, e_{g+i}] = \tau^2$ and with $\tau^{2|k|} = 1$.
  Its abelianization $H_1$ is
  $\bbz^{2g}\oplus\bbz_2$ as desired. When $g$ is even,
  $S_k : B(S,k)\lrar Map_{k-1}(S,S^2)$
  decreases degree by one (see Section \ref{sectionspace})
  but the argument and the conclusion are still the same.
\end{proof}

\bre The above corollary is also a recent calculation of
\cite{bellingo} which is more algebraic in nature and relies on the
full presentation of the braid group $\pi_1(B(S,k))$ for a positive
genus Riemann surface $S$.\ere

\bex We can also apply Proposition \ref{main4} to the case when $M$
is a sphere $S^n$.  Write $\map{}(S^n,S^n) = \coprod_{k\in\bbz}
\map{k}(S^n,S^n)$ for the space of self-maps of $S^n$;
$\map{k}(S^n,S^n)$ being the component of degree $k$ maps.  Since
$\tau^+S^n$ is trivial there is a homeomorphism $\Gamma (\tau^+S^n
)\cong \map{}(S^n,S^n)$. However and as pointed out in \cite{paolo},
one has to pay extra care about components :
$\Gamma_k(\tau^+S^n)\cong \map{k}(S^n,S^n)$ if $n$ is odd and
$\Gamma_k(\tau^+S^n)\cong \map{k-1}(S^n,S^n)$ if $n$ is even (see
section \ref{sectionspace}). Let $p(n)=1$ if $n$ is even and $0$ if
$n$ is odd.  Vassiliev \cite{vassiliev} checks that
$H_*(B(\bbr^n,k);\bbf_2)\lrar H_*(B(\bbr^n,k+1);\bbf_2)$ is an
isomorphism up to degree $k$ and so we get that the map of the
$k$-th braid space of the sphere into the higher free loop space
$$B(S^n,k)\lrar \map{k-p(n)}(S^n,S^n)$$ is a mod-$2$ homology
equivalence up to degree $k-1$. The homology of $\map{}(S^n,S^n)$ is
worked out for all field coefficients in \cite{paolo}. \eex

\bre The braid spaces fit into a filtered construction
$$B(M,n)=: B^1(M,n)\hookrightarrow B^2(M,n)
\hookrightarrow\cdots\hookrightarrow B^n(M,n):=\sp{n}(M)$$ where
$B^p(M,n)$ for $1\leq p\leq n$ is defined to be the subspace
\begin{equation}\label{spnd}
\{[x_1,\ldots, x_n]\in\sp{n}(M)\ |\ \hbox{no more than $p$ of the
$x_i$'s are equal}\}
\end{equation}
Many of our results can be shown to extend with straightforward
changes to $B^p(M,n)$ and $p\geq 1$ when $M$ is a compact Riemann
surface. Some detailed statements and calculations can be found in
\cite{ks}. \ere

\vskip 5pt {\sc Acknowledgements}.  We are grateful to the referee
for his careful reading of this paper.  We would like to thank
Toshitake Kohno, Katsuhiko Kuribayashi and Dai Tamaki for organizing
two most enjoyable conferences first in Tokyo and then in Matsumoto.
Fridolin Roth, Daniel Tanr\'e and Stefan Papadima have motivated
part of this work with relevant questions. We finally thank Fridolin
and Paolo Salvatore for commenting through an early version of this
paper.


\section{Basic Examples and Properties}\label{braids}

As before we write an element of $\sp{n}(X)$ as an unordered
$n$-tuple of points $[x_1,\ldots, x_n]$ or sometimes also as an
abelian finite sum $\sum x_i$ with $x_i\in X$.  For a closed
manifold $M$, $\sp{n}(M)$ is again a manifold for $n>1$ if and only
if $M$ is of dimension less or equal to two \cite{wagner}. We define
$$B(M,n) = \{[x_1,\ldots, x_n]\in\sp{n}(M), x_i\neq x_j, i\neq j\}$$
It is convenient as well to define the "ordered'' $n$-fold
configuration space $F(M,n)= M^n - \Delta_{fat}$ where
\begin{equation}\label{fat}
\Delta_{fat}:= \{(x_1,\ldots, x_n)\in M^n\ |\ x_i=x_j\ \hbox{for
some}\ i=j\}
\end{equation}
is the \textit{fat diagonal} in $M^n$. The configuration space
$B(M,n)$ is obtained as the quotient $F(M,n)/\mathfrak{S}_n$ under
the free permutation action of $\mathfrak{S}_n$ \footnote{In the
early literature on embedding theory \cite{feder}, $B(M,2)$ was
referred to as the ``reduced symmetric square''.}. Both $F(M,n)$ and
$B(M,n)$ are (open) manifolds of dimension $nd$, $d=\dim M$.

Next are some of the simplest non-trivial braid spaces one can
describe.

\ble\label{c2sn} $B(S^n,2)$ is an open $n$-disc bundle over $\bbr
P^n$. When $n=1$, this is the open M\"{o}bius band (see Proposition
\ref{cns1}). \ele

\begin{proof} There is a surjection $\pi : B(S^n,2)\lrar\bbr P^n$ sending
  $[x,y]$ to the unique line $L_{[x,y]}$ passing through the origin and
  parallel to the non-zero vector $x-y$.  The preimage $\pi^{-1}(L_{[x,y]})$
  consists of all pairs $[a,b]$ such that $a-b$ is a multiple of $x-y$.
This can be identified with an ``open'' hemisphere determined by the
hyperplane orthogonal to $L_{[x,y]}$ (i.e. $B(S^n,2)$ can be
identified with the dual tautological bundle over $\bbr P^n$).
\end{proof}

\bex\label{c2rn} Similarly we can see that
$B(\bbr^{n+1},2)\simeq\bbr P^{n}$ and that $B(S^n,2)\hookrightarrow
B(\bbr^{n+1},2)$ is a deformation retract. Alternatively one can see
directly that $B(S^n,2)\simeq\bbr P^n$ for there are an inclusion
$i$ and a retract $r$
\begin{eqnarray*}
i : S^n\hookrightarrow F(S^n,2)\ \ &,&\ \ \  r : F(S^n,2)\lrar
S^n\cr x\longmapsto (x,-x)\ \ &&\ \ \ \ \ \ \  (x,y)\mapsto
{x-y\over |x-y|}
\end{eqnarray*}
Identify $S^n$ with $i(S^n)$ as a subset of $F(S^n,2)$. Then
$F(S^n,2)$ deformation retracts onto this subset via
 $$
 f_t(x,y) = \left({x-ty\over |x-ty|} , {y-tx\over |y-tx|}\right)
 $$
 (which one checks is well-defined). We have that $f_t$ is
 $\bbz_2$-equivariant with respect to the involution $(x,y)\mapsto
 (y,x)$, that $f_0=id$ and that $f_1 : F(S^n,2)\lrar S^n$ is
 $\bbz_2$-equivariant with respect to the antipodal action on $S^n$.  That
 is $S^n$ is a $\bbz_2$-equivariant deformation retraction of $F(S^n,2)$ which
 yields the claim.
\eex

\bex $B(\bbr^2,3)$ is up to homotopy the
  complement of the trefoil knot in $S^3$.  \eex

\bex\label{c2rp2} There is a projection $B(\bbr P^2,2)\lrar\bbr P^2$
which, to any two distinct lines through the origin in $\bbr^3$,
associates the plane they generate and this is an element of the
Grassmann manifold $Gr_2(\bbr^3)\cong Gr_1(\bbr^3) = \bbr P^2$. The
fiber over a given plane parameterizes various choices of two
distinct lines in that plane and that is $B(\bbr P^1,2)=B(S^1,2)$.
As we just discussed, this is an open M\"{o}bius band $M$ and
$B(\bbr P^2,2)$ fibers over $\bbr P^2$ with fiber $M$ (see
\cite{feder}). Interestingly $\pi_1(B(\bbr P^2,2))$ is a quaternion
group of order $16$ (\cite{wang}). \eex

To describe the braid spaces of the circle we can consider the
multiplication map
$$m : \sp{n}(S^1)\lrar S^1\ \ ,\ \
[x_1,\ldots, x_n]\mapsto x_1x_2\cdots x_n$$ Morton \cite{morton}
shows that $m$ is a locally trivial bundle with fiber the closed
$(n-1)$- dimensional disc and this bundle is trivial if $n$ is odd
and non-orientable if $n$ is even. In particular $\sp{2}(S^1)$ is
the closed M\"{o}bius band. In fact one can identify $m^{-1}(1)$
with a closed simplex $\Delta^{n-1}$ so that the configuration space
component $m^{-1}(1)\cap B(S^1,n)$ corresponds to the open part.
This is a non-trivial construction that can be found in
\cite{morton, jack}. Since $B(S^1,n)$ fits in $\sp{n}(S^1)$ as the
open disk bundle one gets that

\bpr\label{cns1} $B(S^1,n)$ is a bundle over $S^1$ with fiber the
open unit disc $D^{n-1}$. This bundle is trivial if and only if $n$
is odd. \epr

Examples \ref{c2rn} and \ref{c2rp2} show that when $\dim M$ is odd
$\neq 1$ or $M$ is not orientable, then $B(M,k)$ fails to be
orientable. The following explains why this needs to be the case.

\ble\label{folklore} (folklore)
  Suppose $M$ is a manifold of dimension $d\geq 2$ and pick $n\geq 2$.  Then
  $B(M,n)$ is orientable if and only if $M$ is orientable of even dimension.
\ele

\begin{proof}
  We consider the $\mathfrak{S}_n$-covering $\pi : F(M,n)\fract{\mathfrak{S}_n}{\lrar}
  B(M,n)$.  If $M$ is not orientable, then so is $M^n$. Now $i :
  F(M,n)\hookrightarrow M^n$ is the inclusion of the complement of
  codimension at least two strata
  so that $\pi_1(F(M,n))\lrar \pi_1(M)^n$ is surjective and hence so is the
  map on $H_1$. The dual map in cohomology is an injection mod $2$ and hence
  $w_1(F(M,n))=i^*(w_1(M^n))\neq 0$ since $w_1(M^n)\neq 0$.  This implies that
  $F(M,n)$ is not orientable if $M$ isn't. It follows that the quotient
  $B(M,n)$ is not orientable as well.

Suppose then that $M$ is orientable. If $d:= \dim M = 2$, then $M$
is a Riemann surface, $B(M,n)$ is open in $\sp{n}(M)$ which is a
complex manifold and hence is orientable. Suppose now that $d:= \dim
M > 2$ so that $\pi_1F(M,n) = \pi_1(M^n)$ (since the fat diagonal
has codimension $>2$). Notice that we have an embedding $\iota :
B(\bbr^d,n)\hookrightarrow B(M,n)$ coming from the embedding of an
open disc $\bbr^d\hookrightarrow M$. Now
$\pi_1(B(\bbr^d,n))=\mathfrak{S}_n$ when $d>2$, and $\iota$ induces
a section of the short exact sequence of fundamental groups for the
$\mathfrak{S}_n$-covering $\pi$ so we have a semi-direct product
decomposition
$$\pi_1(B(M,n)) = \pi_1(M^n) \ltimes\mathfrak{S}_n\ \ \ \ , \ \ \ d > 2$$
Let's argue then that $B(\bbr^d,n)$ is orientable if and only if $d$
is even. Denote by $\tau_x$ the tangent space at $x\in \bbr^d$ and
write $\pi : F(\bbr^d,n)\lrar B(\bbr^d,n)$ the quotient map.  A
transposition $\sigma\in\mathfrak{S}_n$ acts on the tangent space to
$B(\bbr^d,n)$ at some chosen basepoint say $[x_1,\ldots, x_n]$ which
is identified with the tangent space
$\tau_{x_1}\times\cdots\times\tau_{x_n}$ at say $(x_1,\ldots,
x_n)\in\pi^{-1}([x_1,\ldots, x_n])\subset F(\bbr^d,n)\subset
(\bbr^d)^n$.  The action of $\sigma = (ij)$ interchanges both copies
$\tau_{x_i}M$ and $\tau_{x_j}M\cong\bbr^d$ and thus has determinant
$(-1)^d$. Orientation is preserved only when $d$ is even and the
claim follows (for the relation between orientation and fundamental
group see \cite{novikov}, Chapter 4).
\end{proof}

Note that the lemma above is no longer true in the one-dimensional
case according to proposition \ref{cns1}.


\section{Truncated Symmetric Products and Duality}\label{tp}

The heroes here are the truncated symmetric product functors $TP^n$
which were first put to good use in \cite{bcm,lm}. For $n\geq 2$,
define the identification space
$$TP^n(X) := \sp{n}(X)/_{\sim}\ ,\ \
[x,x,y_1\ldots, y_{n-2}]\sim [*,*,y_1,\cdots, y_{n-2}]$$ where as
always $*\in X$ is the basepoint. Clearly $TP^1X = X$ and we set
$TP^0(X) = *$. Note that by adjunction of basepoint $[x_1,\ldots,
x_n]\mapsto [*,x_1,\ldots, x_n]$, we obtain topological embeddings
$\sp{n}(X)\lrar\sp{n+1}X$ and $TP^n(X)\lrar TP^{n+1}X$ of which
limits are $\spy (X)$ and $TP^{\infty}(X)$ respectively. We identify
$\sp{n-1}(X)$ and $TP^{n-1}(X)$ with their images in $\sp{n}(X)$ and
$TP^{n}X$ under these embeddings and we write
\begin{equation}\label{tpnbar}
\overline{TP}^n(X):= TP^n(X)/TP^{n-1}(X)
\end{equation}
for the \textit{reduced} truncated product. These are based spaces
by construction. We will set $\overline{TP}^0(X):= S^0$. The
following two properties are crucial.

\bth\label{property}\
\begin{enumerate}
\item \cite{dt}
$\pi_i(TP^{\infty}(X))\cong\tilde H_i(X;\bbf_2)$
\item \cite{lm} There is a splitting\ \ $H_*(TP^n(X);\bbf_2)\cong
H_*(TP^{n-1}(X);\bbf_2)\oplus \tilde H_*(\overline{TP}^{n}X;\bbf_2 )
$
\end{enumerate}
\end{theorem}

The splitting in (2) is obtained from the long exact sequence for
the pair $(TP^n(X),TP^{n-1}(X))$ and the existence of a retract
$H_*(TP^n(X);\bbf_2)\lrar H_*(TP^{n-1}(X);\bbf_2 )$ constructed
using a transfer argument.  In fact this splitting can be viewed as
a consequence of the following homotopy equivalence discussed in
\cite{lm, zanos}.

\ble $TP^{\infty}(TP^n(X))\simeq TP^{\infty}(\overline{TP}^n(X))
\times TP^{\infty}(TP^{n-1}(X))$. \ele

Further interesting splittings of the sort for a variety of other
functors are investigated in \cite{zanos}. The prototypical and
basic example of course is Steenrod's original splitting of the
homology of symmetric products (which holds with integral
coefficients).

\bth (Steenrod, Nakaoka)\label{steenrod} The induced basepoint
adjunction map on homology \\ $H_*(\sp{n-1}(X);\bbz )\lrar
H_*(\sp{n}(X);\bbz )$ is a split monomorphism.
\end{theorem}


\subsection{Duality and Homological Dimension}
The point of view we adopt here is that $B(M,n) =
TP^n(M)-TP^{n-2}(M)$ as spaces.  A version of Poincar\'e-Lefshetz
duality (Lemma \ref{duality}) can then be used to relate the
cohomology of $B(M,k)$ to the homology of reduced truncated
products. This idea is of course not so new (see \cite{bct, mui}).

If $U\subset X$ is a closed cofibrant subset of $X$, define in
$\sp{n}(X)$ the "ideal"
\begin{equation}\label{ideal}
\underline{U} := \{[x_1,\ldots, x_n]\in\sp{n}(X), x_i\in U\
\hbox{for some $i$}\}
\end{equation}
For example and if $*\in X$ is the basepoint, then $\underline{*} =
\sp{n-1}(X)\subset\sp{n}(X)$.  Let $S$ be the "singular set" in
$\sp{n}(X)$ consisting of unordered tuples with at least two
repeated entries. This is a closed subspace.

\ble\label{quotient} With $U\neq\emptyset$,
$\sp{n}(X)/({\underline{U}\cup S})= \overline{TP}^{n}(X/U)$. \ele

\begin{proof} Denote by $*$ the basepoint of $X/U$ which is the image of $U$
  under the quotient $X\lrar X/U$. Then by inspection
$$\sp{n}(X)/(\underline{U}\cup S) = \sp{n}(X/U)/(\underline{*}\cup S)$$
Moding out $\sp{n}(X/U)$ by $S$ we obtain
  $TP^n(X/U)/TP^{n-2}(X/U)$.  Moding out further
  by $\underline{*}$ we obtain the desired quotient.
\end{proof}

The next lemma is the fundamental observation which states that for
$M$ a compact manifold with boundary and $U\hookrightarrow M$ a
closed cofibration, $B(M-U,n) \cong \sp{n}(M) -
\underline{U\cup\partial M}\cup S$ is Poincar\'e-Lefshetz dual to
the quotient $\sp{n}(M)/(\underline{U\cup\partial M}\cup S)$.  More
precisely, set
\begin{equation}\label{m+}
\overline{M} = M/(U\cup\partial M)
\end{equation}
with the understanding that $\bar M=M$ if $U\cup\partial
M=\emptyset, \{\hbox{point}\}$. The following elaborates on
\cite{bcm}, Theorem 3.2.

\ble\label{duality} If $M$ is a compact manifold of dimension $d\geq
1$, $U\subset M$ a closed subset with $M-U$ connected,
$U\cap\partial M=\emptyset$ and $\overline{M}$ as in (\ref{m+}),
then
$$ H^i(B(M-U,k);\pm\bbz )\cong \begin{cases}
H_{kd-i}(TP^k(\overline{M}), TP^{k-1}(\overline{M});\bbz ),&\
\hbox{if}\ U\cup\partial M\neq\emptyset\\ H_{kd-i}(TP^k(M),
TP^{k-2}(M);\bbz ),&\ \hbox{if}\ U\cup\partial M=\emptyset
\end{cases}
$$
The isomorphism holds with coefficients $\bbf_2$.  When $M$ is even
dimensional and orientable, we can replace $\pm\bbz$ by the trivial
module $\bbz$.  \ele

\begin{proof}
  Suppose $X$ is a compact oriented $d$-manifold with boundary
  $\partial X$.  Then Poincar\'e-Lefshetz duality gives an isomorphism
  $H^{d-q}(X;\bbz )\cong H_q(X,\partial X;\bbz )$.  Apply this to the
  following situation: $X$ is a finite $d$-dimensional CW-complex,
  $V\subset X$ is a closed subset of $X$, and $N$ is a
  tubular neighborhood of $V$ which deformation retracts onto it;
$$V\subset N\subset X$$
$\bar N$ its closure and $\partial\bar N = \partial (X-N)=\bar N-N$.
Assume that $X-N$ is an orientable $d$-dimensional manifold with
boundary $\partial\bar N$. Then
 we have a series of isomorphisms
\begin{equation}\label{iso}
H^{d-q}(X-V;\bbz )\cong H^{d-q}(X-N;\bbz )\cong H_q(X-N,\partial\bar
N;\bbz ) \cong H_q(X,V;\bbz )
\end{equation}

Let's now apply (\ref{iso}) to the case when $X =
\sp{k}(\overline{M})$ with $M$ as in the lemma and with $V$ the
closed subspace consisting of configurations $[x_1,\ldots, x_k]$
such that (i) $x_i=x_j$ for some $i\neq j$ or (ii) for some $i$,
$x_i=*$ the point at which $U\cup\partial M$ is collapsed out.  As
discussed in Lemma \ref{quotient},
$\sp{k}(\overline{M})/\underline{*} =
\sp{k}(M)/(\underline{U\cup\partial M})$ so that
$\sp{k}(\overline{M})/V = \sp{k}(M)/(\underline{U\cup\partial M}\cup
S)$ with $S$ again being the image of the fat diagonal in
$\sp{k}(M)$. Then, according to Lemma \ref{quotient} and to its
proof we see that
$$\sp{k}(\overline{M})/V =
\begin{cases}
TP^k(\overline{M})/TP^{k-1}(\overline{M}),& \hbox
{if $\partial M\neq\emptyset$ or $U\neq\emptyset$} \\
TP^k(M)/TP^{k-2}(M),& \hbox{if $M$ closed and $U=\emptyset$}
\end{cases}
$$
Now $B(M-U,k)\cong \sp{k}(M)- \underline{U\cup\partial M}\cup S
=\sp{k}(\overline{M})-V$ is connected (since $M-U$ is), it is $kd$
dimensional and is orientable if $M$ is even dimensional orientable
(Lemma \ref{folklore}). Applying (\ref{iso}) yields the result in
the orientable case. When $B(M-U,k)$ is non orientable,
Poincar\'e-Lefshetz duality holds with twisted coefficients.
\end{proof}

A version of this lemma has been greatly exploited in \cite{bcm,ks}
to determine the homology of braid spaces and analogs. The following
is immediate.

\bco\label{conR} With $M$, $U\subset M$ as in Lemma \ref{duality},
let
$$R_k = \begin{cases}
conn(TP^k(\overline{M})/TP^{k-1}(\overline{M})), &\hbox{if}\ U\cup\partial M\neq\emptyset,\\
conn(TP^k(M)/TP^{k-2}(M)), &\hbox{if}\ U\cup\partial M=\emptyset
\end{cases}
$$
Then $\cohdim_{\pm\bbz}(B(M-U,k)) = dk-R_k-1$. \eco

Theorem \ref{main3} is now a direct consequence of the following
result.

\ble\label{R} Let $M, U$ and $\overline{M}$ as above,
$r=conn(\overline{M})$ with $r\geq 1$. Then
$$R_k \geq \begin{cases}
k+r-1, &\hbox{if}\ U\cup\partial M\neq\emptyset,\\
k+r-2, &\hbox{if}\ U\cup\partial M=\emptyset
\end{cases}$$
\ele

The proof of this key lemma is based on a computation of Nakaoka
(\cite{nakaoka}, proposition 4.3).  We write $Y^{(k)}$ for the
$k$-fold smash product of a based space $Y$ and $X_{\mathfrak{S}_k}$
the orbit space of a $\mathfrak{S}_k$-space $X$.

\bth\label{nakak} (Nakaoka) If $Y$ is $r$-connected, then
$(Y^{(k)}/\Delta_{fat})_{\mathfrak{S}_k}$ is $r+k-1$ connected.
\end{theorem}

\bre
  In fact nakaoka only proves the homology version of this result and also
  assumes $r\geq 1$. An inspection of his proof shows that $r\geq 0$ works
  as well. Also his homology statement can be upgraded to a genuine
  connectivity statement. To see this, we can assume that $k\geq
  2$ (the case $k=1$ being trivial).  One needs to show in that case that
  $\pi_1((Y^{(k)}/\Delta_{fat})_{\mathfrak{S}_k})=0$.  This follows by an
  immediate application of Van Kampen and the fact that
  $\pi_1(Y^{(k)}/\mathfrak{S}_k)=\pi_1(\bsp{k}Y)=0$ for $k\geq 2$. To see this
  last statement, recall that the natural map $\pi_1(Y)\lrar\pi_1(\sp{k}Y)$
  factors through $H_1(Y;\bbz )$ and then induces an isomorphism
  $H_1(Y;\bbz)\cong\pi_1(\sp{k}Y)$ when $k\geq 2$ \cite{smith}. But if
  $\sp{k-1}(Y)\hookrightarrow \sp{k}(Y)$ induces a surjection on
  fundamental groups, then the cofiber is simply connected (Van-Kampen).
  \ere

\begin{proof} (of Lemma \ref{R} and Theorem \ref{main3}) By construction we
  have the equality $\overline{TP^k}(Y)=
  (Y^{(k)}/\Delta_{fat})_{\mathfrak{S}_k}$.  The connectivity of
  $TP^k(M)/TP^{k-1}(M)$ is (at least) $k+r-1$ according to Theorem
  \ref{nakak}, while that of $TP^{k-1}(M)/TP^{k-2}(M)$ is at least $k+r-2$
  which means that $conn (TP^k(M)/TP^{k-2}(M)) \geq k+r-2$ (by the long exact
  sequence of the triple $(TP^{k-2}(M),TP^{k-1}(M),TP^k(M))$). This produces
  the lower bounds on $R_k$ in Lemma \ref{R}.  Since the cohomology of $B(M-U,k)$
  starts to vanish at $dk-R_k$ (Corollary \ref{conR}), Theorem \ref{main3}
  follows.
\end{proof}


\section{Braid Spaces of Punctured Manifolds}\label{punctured}

We start with a simple proof of Theorem \ref{main}, (\ref{main11});\
$\dim M=d\geq 2$ throughout.

\begin{proof} (Theorem \ref{main}-(\ref{main11}))
This is a direct computation (with $M$ closed)
\begin{eqnarray*}
 H^j(B(M,n);\bbf_2)
&\cong& H_{nd-j}(TP^nM, TP^{n-2}M;\bbf_2 )\ \ \ \ \ \ \ \ \ \ \ \ \
\ \ \ \ \ \ \ \ \ \ \
\hbox{(Lemma \ref{duality})}\\
&\cong& \tilde H_{nd-j}(\overline{TP}^nM;\bbf_2 )\oplus
\tilde H_{nd-j}(\overline{TP}^{n-1}M;\bbf_2 ) \ \ \ \ \ \ (by\ \ref{property}, (2))\\
&\cong&H^j(B(M-\{p\},n);\bbf_2)\oplus H^{j-d}(B(M-\{p\},n-1);\bbf_2
)
\end{eqnarray*}
In this last step we have rewritten $H_{nd-j}$ as $H_{(n-1)d-(j-d)}$
and reapplied Lemma \ref{duality}.
\end{proof}

\noindent{\sc Example}: When $M=S^d$ and $n=2$, then
$B(S^d,2)\simeq\bbr P^d$ and $B(S^d-p,2)=B(\bbr^d,2) = \bbr P^{d-1}$
in full agreement with the splitting. This shows more importantly
that the splitting is not valid for coefficients other than
$\bbf_2$.  The general case is covered by the following observation
of Segal and McDuff.

\ble\label{longexact}\cite{dusa} There is a long exact sequence
$$\rightarrow H_{*-d+1}(B(M-*,n-1))\rightarrow
H_*(B(M-*,n))\rightarrow H_*(B(M,n))\rightarrow
H_{*-d}(B(M-*,n-1))\cdots
$$
\ele

\begin{proof}
  Let $U$ be an open disc in $M$ of radius $<\epsilon$ and let $N =
  M-U$. We have that $B(M-*,n)\simeq B(N,n)$. There is an obvious
  inclusion $B(N,n)\lrar B(M,n)$ and so we are done if we can show
  that the cofiber of this map is $\Sigma^dB(N,n-1)_+$.  To that end
  using a trick as in \cite{dusa} (proof of theorem 1.1) we replace
  $B(M,n)$ by the homotopy equivalent model $B'(M,n)$ of
  configurations $[x_1,\ldots, x_n]\in B(M,n)$ such that at most one
  of the $x_i$'s is in $U$. The cofiber of $B(N,n)\hookrightarrow
  B'(M,n)$ is a based space at $*$ and consists of pairs $(x,D)\in
  \bar U\times B(N,n-1)$ such that if $x\in\partial\bar U$ then
  everything is collapsed out to $*$. But $U\cong D^d$ and $\bar
  U/\partial \bar U = S^d$ so that the cofiber is the half-smash
  product $S^d\rtimes B(N,n-1) = \Sigma^dB(N,n-1)_+$ as asserted.
\end{proof}

In order to prove Theorem \ref{main} we need the following result of
Mostovoy.

\ble\cite{mostovoy}\label{circle} There is a homeomorphism
$\tp{n}(S^1)\cong\bbr P^n$.\ele

\bre We only need that the spaces be homotopy equivalent. It is
actually not hard to see that $\tp{n}(S^1)$ has the same homology as
$\bbr P^n$ since it can be decomposed into cells one for each
dimension less than $n$ and with the right boundary maps. The $k$-th
skeleton is $\tp{k}(S^1)$.  Indeed identify $S^1$ with $[0,1]/\sim$.
A point in $\tp{k}(S^1)$ can be written as a tuple $0\leq
t_1\leq\cdots\leq t_k\leq 1$ with identifications at $t_1=0, t_k=1$
and $t_i=t_{i+1}$. The set of all such points is therefore the image
$\sigma^k$ of a $k$-simplex $\Delta^k\lrar\tp{k}(S^1)$ with
identifications along the faces $F_i\Delta^k$.  Since all faces
corresponding to $t_i=t_{i+1}$ map to the lower skeleton
($\tp{k-2}(S^1))$ and since the last face $F_k\Delta^k$ (when
$t_k=1$) is identified with the zeroth face ($t_1=0$) in
$\tp{k}(S^1)$, the corresponding \textit{chain} map sends the
boundary chain $\partial\sigma_k$ to the image of $\partial\Delta^k
= \sum_{i=0}^k (-1)^iF_i\Delta^k$; that is to the image of
$F_0\Delta^k + (-1)^kF_k\Delta^k$ which is $(1+(-1)^k)\sigma^{k-1}$.
\ere

We need one more lemma.

\ble Set $\overline{TP}^0(X) = S^0$. Then $\overline{TP}^n(X\vee Y)
= \bigvee_{r+s=n} \overline{TP}^r(X) \wedge\overline{TP}^s(Y)$. \ele

\begin{proof}
  Here the smash products are taken with respect to the canonical
  basepoints of the various $\overline{TP}$'s.  A configuration
  $[z_1,\ldots, z_n]$ in $TP^n(X\vee Y)$ can be decomposed into a pair
  of the form $[x_1,\ldots, x_r]\times [y_1,\ldots, y_s]$ in
  $TP^r(X)\times TP^s(Y)$ for some $r+s = n$.  This decomposition is
  unique if we demand that the basepoint (chosen to be the wedgepoint
  $*$) is not contained in the configuration. The ambiguity coming
  from this basepoint is removed when we quotient out $TP^n(X\vee Y)$
  by $\underline{*}=TP^{n-1}(X\vee Y)$, and when we quotient out
  $\bigcup_{r+s = n}TP^r(X)\times TP^s(Y)$ by those pairs of
  configurations with the basepoint in either one of them.  The proof
  follows.
\end{proof}

We are now in a position to prove the second splitting
(\ref{main22})

\begin{proof} (of Theorem \ref{main}-(\ref{main22})) Let $Q_k = \{p_1,\ldots,
  p_k\}$ be a finite subset of $M$ of cardinality $k$. We note that
  the quotient $M/Q_k$ is of the homotopy type of the bouquet
  $M\vee\underbrace{S^1\vee\cdots\vee S^1}_{k-1}$, and that
  $\overline{TP}^l(S^1) = \bbr P^l/\bbr P^{l-1} = S^l$.  Using field
  coefficients we then have (whenever we quote Lemma \ref{duality}
  below we assume that either $M$ is even dimensional orientable or
  that $\bbf=\bbf_2$)
\begin{eqnarray*}
&& H^j(B(M-Q_k,n);\bbf )\\
&\cong& \tilde H_{nd-j}(\overline{TP}^n (M/Q_k))\ \ \ \ \ \ \ \ \ \
\ \ \ \ \ \ \
\hbox{(Lemma \ref{duality} with $U\cup\partial M=Q_k$)}\\
&\cong&
 \tilde H_{nd-j}(\overline{TP}^n(M\vee\bigvee_{k-1}S^1))\\
&\cong& \tilde H_{nd-j}\left(\bigvee_{r+s_1+\cdots + s_{k-1}=n}
\overline{TP}^r(M)\wedge \overline{TP}^{s_1}(S^1)\wedge\cdots\wedge
\overline{TP}^{s_{k-1}}(S^1)\right)\\
&\cong&
 \tilde H_{nd-j}\left(\bigvee_{r+s_1+\cdots + s_{k-1}=n}
S^{n-r}\wedge \overline{TP}^rM \right)\\
&\cong& \bigoplus_{r+s_1+\cdots + s_{k-1}=n}
 \tilde H_{nd-j-n+r}(\overline{TP}^rM )\\
&\cong&\bigoplus_r
 \tilde H_{nd-j-n+r} (\overline{TP}^rM )^{\oplus p(k-1,n-r)}\\
&\cong&\bigoplus_{r=0}^n H^{j - (n-r)(d-1)} (B(M-\{p\},r);\bbf
)^{\oplus p(k-1,n-r)}\ \ \ \ \ \ \ \ \ \hbox{(Lemma \ref{duality})}
\end{eqnarray*}
This is what we wanted to prove.
\end{proof}


\section{Connectivity of Symmetric Products}\label{bounds}

In this section we prove Theorem \ref{connectivity} and Proposition
\ref{conntwo} of the introduction.

\bth\label{connectivity2} Suppose $X$ is a based $r$-connected
  simplicial complex with $r\geq 1$ and let $n\geq 1$. Then
  $\bsp{n}(X)$ is $2n+r-2$- connected.
\end{theorem}

\begin{proof} The claim is tautological for $n=1$ and so we assume throughout
  that $n>1$.  We use some key ideas from \cite{dwyer,kk}.  Start with $X$
  simply connected and choose a CW complex $Y$ such that $H_*(\Sigma
  Y)=H_*(X)$. If $X$ is based and $r$-connected, then $Y$ is based and
  $(r-1)$-connected. A crucial theorem of Dold \cite{dold} now asserts that
  $H_*(\sp{n}X)$, and hence $H_*(\bsp{n}X)$, only depends on $H_*(X)$ so that
  in our case $H_*(\bsp{n}X)=H_*(\bsp{n}\Sigma Y)$. As before we write
  $X^{(n)}$ the $n$-fold smash product of $X$ so that we can identify
  $\bsp{n}X$ with the quotient $X^{(n)}/\mathfrak{S}_n$ by the action of
  $\mathfrak{S}_n$.  It will also be convenient to write
  $X^{(n)}_{\mathfrak{S}_n}:=X^{(n)}/\mathfrak{S}_n$.  Note that $X^{(n)}$ has
  a preferred basepoint which is fixed by the action of $\mathfrak{S}_n$ (i.e
  the action is \textit{based}).  By construction we have equivalences
\begin{equation}\label{bspn}
\bsp{n}(\Sigma Y) = (\Sigma Y)^{(n)}_{\mathfrak{S}_n} = (S^1\wedge
Y)^{(n)}_{\mathfrak{S}_n} =
(S^1)^{(n)}\wedge_{\mathfrak{S}_n}Y^{(n)}
\end{equation}
where here $A\wedge_{\mathfrak{S}_n}B$ is the notation for the
quotient by the diagonal action of $\mathfrak{S}_n$ on $A\wedge B$
where $A$ admits a based right action of $\mathfrak{S}_n$ and $B$ a
based left action.

We next observe that the quotient $(S^1)^{(n)}/K$ is contractible
for any non-trivial subgroup $K\subset\mathfrak{S}_n$. This can be
deduced from the fact (\cite{kk}, section 6) that the permutation
action desuspends in the sense that
$$(S^1)^{(n)}/K\simeq \Sigma \left(S^{n-1}/K\right)$$
where $S^{n-1}$ is viewed as the unit sphere in $\bbr^n$ on which
$K\subset\mathfrak{S}_n$ acts by permutation of coordinates that is
by reflections across the hyperplanes $x_i=x_j$ in $\bbr^n$. Since
$S^{n-1}/K$ is obviously contractible for non-trivial $K$, then so
is $(S^1)^{(n)}/K$. We can then use proposition 7.11 of \cite{dwyer}
to conclude that $(S^1)^{(n)}\wedge_{\mathfrak{S}_n} \Delta_{fat}$
is contractible with $\Delta_{fat}$ as in (\ref{fat}). This subspace
can then be collapsed out in the expression of $\bsp{n}(\Sigma Y)$
of (\ref{bspn}) without changing the homotopy type and one obtains
\begin{equation}\label{finalform}
\bsp{n}(\Sigma Y)\simeq (S^1)^{(n)}\wedge_{\mathfrak{S}_n}
\left(Y^{(n)}/\Delta_{fat}\right)
\end{equation}
The point of expressing $\bsp{n}(X)$ in this form is to take
advantage of the fact that the action of $\mathfrak{S}_n$ on
$Y^{(n)}/\Delta_{fat}$ is based free (i.e. free everywhere but at a
single fixed point say $x_0$ to which the entire $\Delta_{fat}$ is
collapsed out).

Consider the projection $W_n :=
S^{n}\times_{\mathfrak{S}_n}(Y^{(n)}/\Delta_{fat}) \rightarrow
(Y^{(n)}/\Delta_{fat})_{\mathfrak{S}_n}$. This map is a fibration on
the complement of the point $x_0$ with fiber $S^n$ there, and over
$x_0$ the fiber is $F_0 = S^{n}/\mathfrak{S}_n$ (which is
contractible).  The space $\bsp{n}(\Sigma Y)$ in (\ref{finalform})
is obtained from $W_n$ by collapsing out $F_0$ (being contractible
this won't matter) and
$X_n:=*\times_{\mathfrak{S}_n}(Y^{(n)}/\Delta_{fat}) =
(Y^{(n)}/\Delta_{fat})_{\mathfrak{S}_n}$.  Consider the sequence of
maps $(S^n,*)\lrar (W_n,X_n)\fract{}{\lrar} (X_n,X_n)$. This is a
fibration away from the point $x_0\in X$ as we pointed out.  One can
then construct a relative serre spectral sequence (as in \cite{kk},
\S4) with $E^2$-term
$$E^2 = \tilde H_*(X_n; \tilde H_*(S^n))\ \Longrightarrow\
H_*(W_n,X_n)\cong H_*(\bsp{n}(\Sigma Y))$$ But $X_n$ is
$r+n-2$-connected (Theorem \ref{nakak}), $r+n-2\geq 1$, so that the
$E^2$-term is made out of terms of homological dimension $r+n-1+n
=2n+r-1$ or higher which implies that $\bsp{n}(\Sigma Y)=\bsp{n}(X)$
has trivial homology up to $2n+r-2$. But $\bsp{n}(X)$ is simply
connected if $n\geq 2$ (see remark after Theorem \ref{nakak}) and
the proof follows.
\end{proof}

\bex There is a homotopy equivalence \ $\bsp{2}(S^k)\simeq
\Sigma^{k+1}\bbr P^{k-1}$\ (see Hatcher, chapter 4 , example 4K.5).
This space is $k+1 = 4+ (k-1)-2$-connected as predicted and this is
sharp.
 \eex


\subsection{Two dimensional complexes}\label{twodim}
To prove Proposition \ref{conntwo} we use a minimal and explicit
complex constructed in \cite{ks2}. The existence of this complex is
due to the simple but exceptional property in dimension two that
$\sp{n}{D}$, where $D\subset\bbr^2$ is a disc, is again a disc of
dimension $2n$.  Write $X = \bigvee^wS^1\cup (D^2_1\cup\cdots\cup
D^2_r)$ and denote by $\star $ the symmetric product at the chain
level.  In \cite{ks2} we constructed a space $\tsp{n}X$ homotopy
equivalent to $\sp{n}(X)$ and such that $\tsp{}X\simeq
\coprod_{n\geq 0}\tsp{n}X$ has a multiplicative cellular chain
complex generated under $\star $ by a zero dimensional class $v_0$,
degree one classes $e_1,\ldots, e_w$ and degree $2s$ classes
$\sp{s}D_i$, $1\leq i\leq r$, $1\leq s$, under the relations
\begin{eqnarray*}
e_i\star e_j = -e_j\star e_i \,\ (i \neq j)\ \ ,\ \ e_i\star e_i = 0 \ \ \\
\sp{s}D_i\star \sp{t}D_i = {s+t\choose t}\sp{s+t}D_i \,
\end{eqnarray*}
The cellular boundaries on these cells were also explicitly computed
(but we don't need them here). The point however is that a cellular
chain complex for $\tsp{n}(X)$ consists of the subcomplex generated
by cells
$$v_0^r\star e_{i_1}\star \cdots \star e_{i_t}\star
\sp{s_1}(D_{j_1})\star \cdots \star \sp{s_l}(D_{j_l})$$ with
$r+t+s_1+\cdots +s_l= n$ and $t\leq w$ where $w$ again is the number
of leaves in the bouquet of circles. The dimension of such a cell is
$t+2(s_1+\cdots + s_l)$ for pairwise distinct indices among the
$e_i$'s.

A cellular complex for $\bsp{n}X$ can then be taken to be the
quotient of $C_*(\tsp{n}X)$ by the summand $C_*(\tsp{n-1}X)$ and
this has cells of the form
$$e_{i_1}\star \cdots \star e_{i_t}\star \sp{s_1}(D_{j_1})\star \cdots \star \sp{s_l}(D_{j_l})$$
with $t+s_1+\cdots +s_l= n$. The dimension of such a cell is
$t+2(s_1+\cdots + s_l)=2n-t$.  The smallest such dimension is
$2n-min(w,n)$. This means that $conn(\tsp{n}X/\tsp{n-1}X) =
conn(\bsp{n}X) \geq 2n-min(w,n)-1$ and Proposition \ref{conntwo}
follows.

\bex
  A good example to illustrate Proposition \ref{conntwo} is when $S$ is a
  closed Riemann surface of genus $g$. It is well-known that for $n\geq 2g-1$,
  $\sp{n}(S)$ is an analytic fiber bundle over the Jacobian (by a result of
  Mattuck)
$$\bbp^{n-g}\lrar\sp{n}(S)\fract{\mu}{\lrar} J(S)$$
where $\mu$ is the Abel-Jacobi map. In fact this is the
projectivisation of an $n-g+1$ complex vector bundle over $J(S)$.
Collapsing out fiberwise the hyperplanes
$\bbp^{n-g-1}\subset\bbp^{n-g}$ we get a fibration $\zeta_n :
S^{2n-2g}\lrar E_n\lrar J(S)$ with a preferred section, so that for
$n\geq 2g$, $\bsp{n}(S)$ is the cofiber of this section. This is
$2n-2g-1$-connected as predicted, and in fact $\tilde
H_*(\bsp{n}(S))= \sigma^{2n-2g}H_*(J(S))$ where $\sigma$ is a formal
suspension operator which raises degree by one. \eex

\subsection{Connectivity and Truncated Products}\label{spec}
The homology of truncated products, and hence of braid spaces, is
related to the homology of symmetric products via a very useful
spectral sequence introduced in \cite{bcm}. This spectral sequence
has been used and adapted with relative success to other situations;
eg. \cite{ks}. The starting point is the duality in lemma
\ref{duality}. The problem of computing $H^*(B(M,n);\bbf )$ becomes
then one of computing the homology of the relative groups
$H_{*}(TP^n\overline{M}, TP^{n-2}\overline{M};\bbf )$. The key tool
is the following \textit{Eilenberg-Moore type} spectral sequence
with field coefficients $\bbf$.

\bth\label{specseq}\cite{bcm}
  Let $X$ be a connected space with a non-degenerate basepoint. Then there is
  a spectral sequence converging to $H_{*}(TP^{n}(X), TP^{n-1}(X);\bbf )$ ,
  with $E^1$-term
\begin{equation}\label{unpunctured}
  \bigoplus_{i+ 2j= n}
  H_*(\sp{i}X,\sp{i-1}X)\tensor H_*(\sp{j}(\Sigma X),\sp{j-1}(\Sigma X))
\end{equation}
and explicit $d^1$ differentials.
\end{theorem}

Field coefficients are used here because this spectral sequence uses
the Kunneth formula to express $E^1$ as in (\ref{unpunctured}). Here
$\sp{-1}(X)=\emptyset$ and $\sp{0}(X)$ is the basepoint.

\bex When $X=S^1$, then $H_*(TP^{n}(S^1),TP^{n-1}(S^1))=\tilde
  H_*(S^n)$.  Since $\sp{i}S^1\simeq S^1$ for all $i\geq 1$, the spectral
  sequence in this case has $E^1$-term of the form
$$H_*(S^1,*)\tensor H_*(\bbp^{{n-1\over 2}},\bbp^{{n-1\over 2}-1})
= \sigma\tilde H_*(S^{n-1}) = \tilde H_*(S^n)$$ if $n$ is odd (where
$\sigma$ is the suspension operator), or $E^1_{*,*}=
H_*(\bbp^{(n/2)},\bbp^{(n/2)-1})= \tilde H_*(S^n)$ if $n$ is even.
In all cases the spectral sequence collapses at $E^1$. \eex

Now Lemma \ref{duality} combined with Theorem \ref{specseq} gives an
easy method to produce upper bounds for the non-vanishing degrees of
$H^*(B(M,n))$. The least connectivity of the terms
$\bsp{i}X\times\bsp{j}(\Sigma X)$ for $i+2j=n$ translates by duality
to such an upper bound.  This was in fact originally our approach to
the cohomological dimension of braid spaces.  We illustrate how we
can apply this spectral sequence by deriving Corollary
\ref{twocomplexes} from Proposition \ref{conntwo}.

\begin{proof} (Corollary \ref{twocomplexes})
Suppose $Q\cup\partial S\neq\emptyset$. The spectral sequence of
Theorem \ref{specseq} converging to the homology of
$(TP^k(\overline{S}),TP^{k-1}(\overline{S}))$ takes the form
\begin{equation}\label{e1}
  E^1 = \tilde H_*(\bsp{k}\overline{S})
  \bigoplus\oplus_{i+2j=k}(H_*(\bsp{i}\overline{S})\otimes
  H_*(\bsp{j}(\Sigma \overline{S}))
  \bigoplus \tilde H_*(\bsp{k/2}(\Sigma \overline{S}))
\end{equation}
(if $k$ odd, the far right term is not there). We have that $R_k$
(as in Corollary \ref{conR}) is at least the connectivity of this
$E^1$-term. Since $\overline{S}$ is a two dimensional complex, the
connectivity of $\bsp{i}(\overline{S})$ is at least $2i-min(w,i)-1$
(for some $w\geq 0$). The connectivity of $\bsp{j}(\Sigma
\overline{S})$ is at least $2j + r-2\geq 2j-1$ since $\Sigma
\overline{S}$ is now simply connected (Theorem \ref{connectivity2}).
The connectivity of $\bsp{i}(\overline{S})\wedge \bsp{j}(\Sigma
\overline{S})$ for non-zero $i$ and $j$ is then at least
$$(2i-min(w,i)-1)+(2j-1)+1 = i+k-min(w,i)-1$$
When $i=0$, then $j={k\over 2}$ ($k$ even) and
$conn(\bsp{k/2}(\Sigma \overline{S}))\geq k-1$. The connectivity of
the $E^1$-term (\ref{e1}) is at least the minimum of
$$
\begin{cases}
i+k-min(w,i)-1  &  1\leq i\leq k-1\\
2k-min(w,k)-1& i=k\\
k-1& i=0
\end{cases}
$$
which is $k-1$. By duality $H^*(B(S-Q,k)) = 0$ for $* \geq 2k-k+1 =
k+1$. If $S$ is closed, then the same argument shows that this bound
needs to be raised by one.
\end{proof}


\section{Stability and Section Spaces}\label{stability}

In this final section, we extrapolate on standard material and make
slightly more precise a well-known relationship between
configuration spaces and section spaces
\cite{dusa,bct,segal1,quarterly}.

When manifolds have a boundary or an end (eg. a puncture), one can
construct embeddings
\begin{equation}\label{marching}
+ : B(M,k)\lrar B(M,k+1)
\end{equation}
by "addition of points" near the boundary, near "infinity" or near
the puncture. In the case when $\partial M\neq\emptyset$ for
example, one can pick a component $A$ of the boundary and construct
a nested sequence of collared neighborhoods $V_1\supset
V_2\supset\cdots \supset A$ together with sequences of points
$x_k\in V_{k}-V_{k+1}$. There are then embeddings $B(M-V_k,k)\lrar
B(M-V_{k+1},k+1)$ sending $\sum z_i$ to $\sum z_i+x_k$. Now we can
replace $B(M-V_k,k)$ by $B(M-A,k)$ and then by $B(M,k)$ up to small
homotopy.  In the direct limit of these embeddings we obtain a space
denoted by $B(M,\infty )$.  Note that an easy analog of Steenrod's
splitting \cite{bcm} gives the splitting
\begin{equation}\label{split2}
H_*(B(M,\infty ))\cong\bigoplus_{k=0} H_*(B(M,k+1), B(M,k))
\end{equation}
(here $B(M,0)=\emptyset$).  In fact (\ref{split2}) is a special case
of a trademark \textit{stable splitting} result for configuration
spaces of open manifolds or manifolds with boundary.  Denote by
$D_k(M)$ the cofiber of (\ref{marching}). For example
$D_1(M)=B(M,1)=M$.

\bth\label{split3} \cite{bodig,cohen} For $M$ a manifold with
non-empty boundary, there is a stable splitting (i.e. after
sufficiently many suspensions)
$$B(M,k)\simeq_s\bigvee_{i=0}^kD_i(M)$$
\end{theorem}

The classical case of $M=D^n$ (closed $n$-ball) is due to Victor
Snaith.  A short and clever argument of proof for this sort of
splittings is due to Fred \cite{cohen}. The next stability bound is
due to Arnold and a detailed proof is in an appendix of
\cite{segal1}.

\bth \label{arnold} (Arnold) The embedding $B(M,k)\hookrightarrow
  B(M,k+1)$ induces a homology monomorphism and a homology equivalence up to
  degree $[k/2]$.
\end{theorem}

The monomorphism statement is in fact a consequence of
(\ref{split2}). Arnold's range is not optimal. For instance

\bth\cite{ks} If $S$ is a compact Riemann surface and $S^*=S-\{p\}$,
then $B(S^*,k)\hookrightarrow B(S^*,k+1)$ is a homology equivalence
up to degree $k-1$.
\end{theorem}

We define $s(k)$ to be the homological connectivity of $+ :
B(M,k)\lrar B(M,k+1)$ (see \S1.3) .  By Arnold, $s(k)\geq [k/2]$ .


\subsection{Section Spaces}\label{sectionspace}
If $\zeta : E\lrar B$ is a fiber bundle over a base space $B$, we
write $\Gamma (\zeta )$ for its space of sections.  If $\zeta$ is
trivial then evidently $\Gamma(\zeta )$ is the same as maps into the
fiber.  Let $M$ be a closed smooth manifold of dimension $d$,
$U\subset M$ a closed subspace and $\tau^+M$ the fiberwise one-point
compactification of the tangent bundle over $M$ with fiber $S^d =
\bbr^d \cup\{\infty\}$. Then $\tau^+M\lrar M$ has a preferred
section $s_{\infty}$ which is the section at $\infty$ and we let
$\Gamma(\tau^+M;U )$ be those sections which coincide with
$s_{\infty}$ on $U$. Note that $\Gamma (\tau^+M )$ splits into
components indexed by the integers as in
$$\Gamma (\tau^+M) := \coprod_{k\in\bbz} \Gamma_k(\tau^+M)$$
This degree arises as follows. Let $s : M\lrar\tau^+M$ be a section.
By general position argument it intersects $s_{\infty}$ at a finite
number of points and there is a sign associated to each point.  This
sign is defined whether the manifold is oriented or not (as in the
definition of the Euler number). The degree is then the signed sum.
Similarly we can define a (relative) degree of sections in $\Gamma
(\tau^+M;U)$.

Observe that if $\tau^+M$ is trivial, then $\Phi : \Gamma
(\tau^+M)\fract{\simeq}{\lrar}\map{}(M,S^d)$, where $d=\dim M$. The
components of $\map{}(M,S^d)$ are indexed by the degree of maps
(Hopf), but at the level of components we have the equivalence
$$\Gamma_{k} (\tau^+M)\simeq\map{k+\ell}(M,S^d)$$ where $\ell$ is such
that $\Phi (s_{\infty})\in\map{\ell}$.  In the case when
$M=S^{even}$, then $\Phi (s_{\infty})$ is the antipodal map which
has degree $\ell = -1$ \cite{paolo}. When $M=S$ is a compact Riemann
surface, $\ell =-1$ when the genus is even and $\ell = 0$ when the
genus is odd \cite{ks}. Further relevant homotopy theoretic
properties of section spaces are summarized in the appendix.

\subsection{Scanning and Stability}\label{scan}
A beautiful and important connection between braid spaces and
section spaces can be found for example in
\cite{segal2,dusa,quarterly} (see \cite{james} for the fiberwise
version).  This connection is embodied in the "scanning" map
\begin{equation}\label{scanning}
S_k : B(M-U,k)\lrar \Gamma_k(\tau^+M; U\cup\partial M )
\end{equation}
where $U$ is a closed subspace of $M$. Here and throughout we assume
that removing a subspace as in $M-U$ doesn't disconnect the space.
The scanning map has very useful homological properties. A sketch of
the construction of $S_k$ for closed Riemaniann $M$ goes as follows
(for a construction that works for topological manifolds see for
example \cite{dww}). First construct $S_1: M-U\lrar\Gamma_1$. We can
suppose that $M$ has a Riemannian metric and use the existence of an
exponential map for $\tau M$ which is a continuous family of
embeddings $exp_x: \tau_xM\lrar M$ for $x\in M$ such that $x\in
im(exp_x)$ and $im(exp_x)^+\cong\tau_x^+M$ (the fiber at $x$ of
$\tau^+M$).  By collapsing out for each $x$ the complement of
$im(exp_x)$ we get a map $c_x : M\lrar im(exp_x)^+\cong\tau_x^+M$
Let $V$ be an open neighborhood of $U$, $M-V\lrar M-U$ being a
deformation retract. Then we have the map
$$S_1 : M-V\lrar\Gamma (\tau^+M)\ ,\ y\mapsto (x\mapsto c_x(y))\in
\tau_x^+M$$ Observe that for $x$ near $U$, the section $S_1(y)$
agrees with the section at infinity (i.e. we say it is
\textit{null}). In fact and more precisely, $S_1$ maps into
$\Gamma^c(\tau^+M,U)$ the space of sections which are null outside a
compact subspace of $M-U$.  A deformation argument shows that
$\Gamma^c\simeq\Gamma$.  It will be convenient to say that a section
$s\in\Gamma$ is \textit{supported} in a subset $N\subset M$ if
$s=s_{\infty}$ outside of $N$. A useful observation is that if
$s_1,s_2$ are two sections supported in closed $A$ and $B$ and
$A\cap B=\emptyset$, then we can define a new section which is
supported in $A\cup B$, restricting to $s_1$ on $A$ and to $s_2$ on
$B$.

Extending $S_1$ to $S_k$ is now easy.  We first choose $\epsilon >
0$ so that $B^{\epsilon}(M,k)$ the closed subset of $B(M,k)$ where
particles have pairwise separation $\geq 2\epsilon$ is homotopic to
$B(M,k)$ (this is verified in \cite{dusa}, lemma 2.3). We next
choose the exponential maps to be supported in neighborhoods of
radius $\epsilon$.  Given a finite subset $Q:=\{y_1,\ldots, y_k\}\in
B^{\epsilon}(M-U,k)$, each point $y_i$ determines a section
supported in $V_i := im(exp_{y_i} )$.  Since the $V_i$'s are
pairwise disjoint, these sections fit together to give a section
$s_Q$ supported in $\bigcup V_i$ so that $S_k(Q):= s_Q$.

When $M$ is compact with boundary, then we get the map in
(\ref{scanning}) by replacing $B(M-U,k)$ by $B(M-U\cup\partial M,k)$
and $\Gamma^c(\tau^+M,U)$ by $\Gamma (\tau^+M, U\cup\partial M)$ the
space of sections that are null outside a compact subspace of
$M-U\cup\partial M$. We let $s(k)$ be the stability range of the map
$B(M-U,k)\lrar B(M-U,k+1)$ (as in \S6.1)

The next proposition is a follow up on a main result of \cite{dusa}
(see also \cite{quarterly}).

\bpr\label{dusa1} Suppose $M$ is a closed manifold and $U\subset M$
a non-empty closed subset, $M-U$ connected. Then the map $S_{k*}:
H_*(B(M-U,k))\lrar H_*(\Gamma_k(\tau^+M,U))$ is a monomorphism in
all dimensions and an isomorphism up to dimension $s(k)$. \epr

\begin{proof} It is easy to see that the maps $S_k$ for various $k$ are
  compatible up to homotopy with stabilization so we obtain a map $S :
  B(M,\infty )\lrar \Gamma_{\infty}(\tau^+M,U):=\lim_k\Gamma_k(\tau^+M,U)$
which according to the main
  theorem of McDuff is a homology equivalence (in fact all components of
  $\Gamma (\tau^+M,U)$ are equivalent and $\Gamma_{\infty}$ can be chosen to
  be the component containing $s_{\infty}$).  But according to (\ref{split2})
  $H_*(B(M-U,k))\rightarrow H_*(B(M-U,\infty ))$ is a monomorphism, and then an
  isomorphism up to dimension $s(k)$. The claim follows.
\end{proof}

This now also implies our last main result from the introduction.

\begin{proof} (of Proposition \ref{main4}) Suppose that $M$ is a closed
manifold of dimension $d$, $U$ a small open neighborhood of the
basepoint $*$
 and consider the fibration (see appendix)
$$\Gamma_k (\tau^+M;\bar U)\lrar\Gamma_k(\tau^+M)\lrar S^d$$
The main point is to use the fact as in (\cite{dusa}, proof of
theorem 1.1) that scanning sends the exact sequence in Lemma
\ref{longexact} to the Wang sequence of this fibration.  Let $N=M-U$
so that we can identify $\Gamma_k (\tau^+M;\bar U)$ with
$\Gamma_k(\tau^+N;\partial N)$ which we write for simplicity
$\Gamma^c_k(\tau^+N)$ as before. Under these identifications and by
a routine check we see that scanning induces commutative diagrams
$$\small
\begin{matrix}
\rightarrow&H_{q-d+1}(B(N,k-1))&\rightarrow&
H_q(B(N,k))&\rightarrow&H_q(B(M,k))&\rightarrow&H_{q-d}(B(N,k-1))
&\rightarrow\\
&\decdnar{S}&&\decdnar{S}&&\decdnar{S}&&\decdnar{S}&\\
\rightarrow&H_{q-d+1}(\Gamma^c_{k}(\tau^+N))&\rightarrow&
H_q(\Gamma^c_{k}(\tau^+N))&\rightarrow&H_q(\Gamma_{k}(\tau^+M))
&\rightarrow&H_{q-d}(\Gamma^c_{k}(\tau^+N)) &\rightarrow
\end{matrix}
$$
where the top sequence is the homology exact sequence for the pair
$(B(M,k),B(N,k))$ as discussed in Lemma \ref{longexact} and the
lower exact sequence is the Wang sequence of the fibration
$\Gamma_k(\tau^+M)\lrar S^d$. According to Proposition \ref{dusa1},
the map $S_{k*}: H_q(B(N,k))\lrar H_q(\Gamma^c_k(\tau^+N))$ is an
isomorphism up to degree $q=s(k)$. It follows that all vertical maps
in the diagram above involving the subspace $N$ together with the
next map on the right (which doesn't appear in the diagram) are
isomorphisms whenever $q\leq s(k-1)\leq s(k)$. By the $5$-lemma the
middle map is then an isomorphism within that range as well. This
proves the proposition.
\end{proof}

We can say a little more when $k=1$, $M$ closed always.

\ble\label{s1} The map $S_1 :M\lrar \Gamma_{1} (\tau^+M)$ induces a
monomorphism in homology in degrees $r+1,r+2$, where $r=conn(M)$,
$r\geq 1$. \ele

\begin{proof}
Consider $\Gamma(s\tau^+M)$ the space of sections of the fibration
$s\tau^+M\lrar M$ obtained from $\tau^+M$ by applying fiberwise the
functor $\spy$. It is easy to see that scanning has a stable analog
$st: \spy (M_+)\lrar \Gamma (s\tau^+M)$ but harder to verify that
$st$ is a (weak) homotopy equivalence \cite{dww, quarterly}. Note
that $\spy (M_+)\simeq\spy M\times\bbz$ and $\spy (M)$ is equivalent
to a connected component (any of them) say $\Gamma_0 (s\tau^+M )$.
By construction the following diagram homotopy commutes
$$
\begin{matrix}\label{fromstosp}
M&\fract{S_1}{\lrar}&\Gamma_{1} (\tau^+M)\\
\decdnar{}&&\decdnar{\alpha}\\
\spy (M)&\fract{st}{\lrar}&\Gamma_0 (s\tau^+M )
\end{matrix}
$$
where the right vertical map $\alpha$ is induced from the natural
fiber inclusion $\alpha : S^d\hookrightarrow\spy (S^d)$. When $M$ is
$r$-connected, the map $M\lrar \spy (M)$ induces an isomorphism in
homology in dimensions $r+1$ and $r+2$ (\cite{nakaoka}, Corollary
4.7).  This means that the composite $M\rightarrow \Gamma_{1}
(\tau^+M)\rightarrow \Gamma_1 (s\tau^+M )$ is a homology isomorphism
in those dimensions and the claim follows.
\end{proof}

\bre If $M$ has boundary, then by scanning $M_0:=M-\partial M$ we
  obtain a map into the compactly supported sections
$\Gamma (\tau^+M)$.  This map extends to a map $S :
  M/\partial M\lrar \Gamma (\tau^+M)$ which is
  according to \cite{aouina} $(d-r+1)$-connected if $M$ is $r$-connected of
  dimension $d\geq 2$.
\ere

\section{Appendix: Some Homotopy Properties of Section Spaces}

All spaces below are assumed connected. We discuss some pertinent
statements from \cite{switzer}. Let $p: E\lrar B$ be a Serre
fibration, $i: A\hookrightarrow X$ a cofibration ($A$ can be empty)
and  $u: X\lrar E$ a given map. Slightly changing the notation in
that paper, we define
$$\Gamma_u (X,A; E,B) = \{f: X\lrar E\ |\ f\circ i = u\circ i,
p\circ f = p\circ u\}$$ This is a closed subspace of the space of
all maps $\map{}(X,E)$ and is in other words the solution space for
the extension problem
$$\xymatrix{
A\ar[r]^{ui}\ar[d]^i&E\ar[d]^p\\
X\ar[r]_{pu}\ar[ru]^u&B }
$$
with data $u_{|A}: A\lrar E$ and $pu : X\lrar B$. When $A=\{x_0\}$
and $B = \{y_0\}$ then $\Gamma (X,x_0;E,y_0) = \bmap{}(X,E)$ is the
space of based maps from $X$ to $Y$ sending $x_0$ to $y_0$.  On the
other hand and when $X=B$ and $A=\emptyset$, then
$\Gamma_u(B,\emptyset;E,B) = \Gamma (E)$ is the section space of the
fibration $\zeta = (E\fract{p}{\lrar} B)$.

\bpr\cite{switzer}\label{switzer}
\begin{itemize}
\item
If $A\subset X'\subset X$ is a nested sequence of NDR pairs, and $j:
X'\hookrightarrow X$ the inclusion, then the induced map
$\Gamma_u(X,A; E,B)\lrar \Gamma_{uj}(X',A;E,B)$\ yields a fibration
with $\Gamma_u (X,X';E,B)$ as fibre.
\item
If $E\lrar E'\lrar B$ are two fibrations and $q: E\lrar E'$ the
projection, then the induced map $\Gamma_u(X,A; E,B)\lrar
\Gamma_{qu}(X,A;E',B)$\ is a fibration with $\Gamma_u (X,A;E,E')$ as
fibre.
\end{itemize}
\epr

The first part of Switzer's result implies that restriction of the
bundle $\zeta : E\lrar B$ to $X\subset B$ is a fibration $\Gamma
(\zeta )\lrar \Gamma (\zeta_{|X})$ with fiber the section space
$\Gamma (\zeta, X)$ i.e. those sections of $\zeta$ which are
"stationary" over $X$ (compare \cite{james}, chapter 1, \S8).  An
example of relevance is when $\zeta = \tau^+M$ is the fiberwise
one-point compactification and $s_{\infty}$ is the section mapping
at infinity.  Denote by $S^d$ the fiber over $x_0\in M$. If $U$ is a
small open neighborhood of $x_0$, then $\Gamma (\zeta_{|\bar
U})\simeq S^d$ and we have a fibration
\begin{equation}\label{fibration}
\Gamma (\tau^+M,\bar U)\lrar\Gamma (\tau^+M)\fract{res}{\lrar} S^d
\end{equation}
where the fiber consists of those sections which coincide with
$s_{\infty}$ on $U$. So for instance if $M=S^d$, $\Gamma
(\tau^+M,\bar U)\simeq\Omega^dS^d$ and the fibration reduces to the
evaluation fibration $\Omega^dS^d\rightarrow
\map{}(S^d,S^d)\rightarrow S^d$.

Finally and according to \cite{james} (p. 29), if $E\lrar B$ is a
Hurewicz fibration and $s,t$ are two sections, then $s$ and $t$ are
homotopic if and only if they are section homotopic. We use this to
deduce the following lemma.

\ble
  \label{mono} Let $\pi : E\lrar B$ be a fibration with a preferred section
  $s_{\infty}$ (which we choose as basepoint). Then the inclusion $\Gamma
  (E)\lrar\map{}(B,E)$ induces a monomorphism on homotopy groups.
\ele

\begin{proof} We give $\Gamma (E)\subset\map{}(B,E)$ the common basepoint
  $s_{\infty}$.  An element of $\pi_i\Gamma (E)$ is the homotopy class of a
  (based) map $\phi : S^i\lrar\Gamma (E)$ or equivalently a map $\phi :
  S^i\times B\lrar E$ (where $\phi (-,b)\in\pi^{-1}(b)$ and $\phi
  (N,-)=s_\infty(-)$, $N$ the north pole of $S^i$) and the homotopy is through
  similar maps. Write $\Phi$ the image of $\phi$ via the composite
  $S^i\lrar\Gamma (E)\lrar\map{}(B,E)$. Now $\Phi$ can be viewed as a section
  of $S^i\times E\lrar S^i\times B$ and a null-homotopy of
  $\Phi$ is a homotopy to $id\times s_{\infty}$. Since this null-homotopy can
  be done fiberwise it is a null-homotopy in $\Gamma (E)$ from $\phi$ to
  $s_{\infty}$.
\end{proof}


\addcontentsline{toc}{section}{Bibliography}
\bibliography{biblio}

\begin{thebibliography}{doC}

\bibitem{aouina} M. Aouina, J. Klein, \textit{On C.T.C. Wall's suspension
    theorem}, On C. T. C. Wall's suspension theorem.  Forum Math.  18 (2006),
  no. 5, 829--837.
\bibitem{dwyer} G. Arone, B. Dwyer, \textit{Partition complexes, Tits
    buildings and symmetric products}, Proc. London Math. Soc. (3) 82 (2001),
  no. 1, 229--256.
\bibitem{atiyah} M. Atiyah, J. Jones, \textit{Topological aspects of
    Yang-Mills theory}, Commun. Math. Phys. \textbf{61} (1978), 97--118.
\bibitem{bellingo} P. Bellingeri, S. Gervais, J. Guaschi, \textit{Lower
    central series for surface braid groups}, preprint Gennaio 2006,
  universita di pisa.
\bibitem{bodig}  C.F. Bodigheimer, \textit{Stable splittings of mapping
    spaces}, Algebraic topology, Proc. Seattle (1985), Springer lecture notes
  \textbf{1286}, 174--187.
\bibitem{bcm}  C.F. Bodigheimer, F.R. Cohen, R.J. Milgram,
  \textit{Truncated symmetric products and configuration spaces}, Math. Zeit.,
  \textbf{214} (1993), 179--216.
\bibitem{bct}  C.F. Bodigheimer, F.R. Cohen, L. Taylor, \textit{On the
    homology of configuration spaces}, Topology, \textbf{48} (1989), 111--123.
\bibitem{cohen} F. Cohen, \textit{The unstable decomposition of
    $\Omega^2\Sigma^2X$ and its applications}, Math. Z.  \textbf{182} (1983),
  no. 4, 553--568.  `\bibitem{cohen2} F. Cohen, \textit{On the mapping class
    groups for punctured spheres, the hyperelliptic mapping class groups,
    ${\rm SO}(3)$, and ${\rm Spin}\sp c(3)$}, Amer. J. Math. {\bf 115} (1993),
  no. 2, 389--434.
\bibitem{james} M. Crabb, I. James, \textit{Fibrewise homotopy theory}, Spring
  Monographs in Math. (1998).
\bibitem{dold} A, Dold, \textit{Homology of symmetric products and other
    functors of complexes}, Ann. of Math. (2) 68 1958 54--80.
\bibitem{dt} A. Dold, R. Thom, \textit{Quasifaserungen and unendliche
    symmetrische produkte}, Ann. Math. {\bf 67}(1958), 239--2
\bibitem{novikov} B. Dubrovin, A. Fomenko, S. Novikov, \textit{Modern
    geometry: methods and applications. Part II}, Springer GTM
  \textbf{93} (1992).
\bibitem{dww} W. Dwyer, M. Weiss, B. Williams, \textit{A parameterized index
    theorem for the algebraic $K$-theory Euler class}, Acta Math.
  \textbf{190} (2003), no. 1, 1--104.
\bibitem{feder} S. Feder, \textit{The reduced symmetric product of projective
    spaces and the generalized Whitney theorem}, Illinois J. Math.
  \textbf{16} (1972), 323--329.
\bibitem{goryunov} V.V. Gorjunov, \textit{Cohomology of braid groups of the
    series C and D}, Trans. Moscow. Math. Soc. {\bf 2} (1982), 233--241.
\bibitem{contemp} S. Kallel, \textit{Configuration spaces and the topology of
    curves in projective space}, Contemp. Math. {\bf 279} (2001), 135--149.
\bibitem{quarterly} S. Kallel, \textit{Particle spaces on manifolds and
    generalized Poincar\'e dualities}, Quart. J. of Math \textbf{52} (2001),
  45--70.
\bibitem{kk} S. Kallel, R. Karoui, \textit{Symmetric joins and weighted
    barycenters}, math.AT/0602283.
\bibitem{ks2} S. Kallel, P. Salvatore, \textit{Symmetric products of two
    dimensional complexes}, Contemp. Math. \textbf{407} (2006), 147--161.
\bibitem{ks} S. Kallel, P. Salvatore, work in progress on projective sigma
  models.
\bibitem{lm} P. Loeffler, R.J. Milgram, \textit{The structure of deleted
    symmetric products}, Braids (Santa Cruz 1986), Contemp. Math.  \textbf{
    78} (1988), 415--424.
\bibitem{dusa} D. McDuff, \textit{Configuration spaces of positive and
    negative particles}, Topology \textbf{ 14}(1975), 91--107.
\bibitem{jack} J. Morava, \textit{The tangent bundle of an almost complex free
    loopspace}, Homology Homotopy Appl. \textbf{3} (2001), no. 2, 407--415.
\bibitem{morton} H. R. Morton, \textit{Symmetric products of the circle},
  Proc.  Camb. Phil. Soc. \textbf{ 63} (1967), 349--352.
\bibitem{mostovoy} J. Mostovoy, \textit{Geometry of truncated symmetric
    products and real roots of real polynomials}, Bull. London Math. Soc.
  \textbf{ 30} (1998), no. 2, 159--165.
\bibitem{mui} H. Mui, \textit{Duality in the infinite symmetric products},
Acta Mathematica Vietnamica, \textbf{5}, no.1 (1980), 100--149.
\bibitem{nakaoka} M. Nakaoka, \textit{Cohomology of symmetric products}, J.
  Institute of Polytechnics, Osaka city University \textbf{8}, no. 2,
  121--140.
\bibitem{nap1} F. Napolitano, \textit{Configuration spaces on surfaces}, C.R.
  Acad. Sci. Paris \textbf{327} (1998) no. 10, 887--892.
\bibitem{nap2} F. Napolitano, \textit{On the cohomology of configuration
    spaces on surfaces}, J. London Math. Soc. (2) \textbf{ 68} (2003),
  477--492.
\bibitem{ossa} E. Ossa, \textit{On the cohomology of configuration spaces},
  Algebraic topology: new trends in localization and periodicity, 353--361,
  Progr. Math., \textbf{136}, Birkhauser 1996.
\bibitem{frido} F. Roth, \textit{On the category of euclidean configuration
spaces}, these proceedings.
\bibitem{paolo} P. Salvatore, \textit{Configuration spaces and higher free
    loop spaces}, Math. Z. \textbf{248} (2004), 527--540.
\bibitem{segal1} G. Segal, \textit{The topology of spaces of rational
    functions}, Acta. Math., \textbf{143}(1979), 39--72.
\bibitem{segal2} G. Segal, \textit{configuration spaces and iterated loop
    spaces}, Inv. \textbf{21} (1973), 213--221.
\bibitem{smith} Smith P.A , \textit{Manifolds with abelian fundamental
groups}, Annals of Mathematics, {\bf 37} (1936), 526--533.
\bibitem{switzer}  R. Switzer, \textit{Counting elements in homotopy
    sets}, Math. Z. \textbf{178} (1981), 527--554.
\bibitem{vassiliev} V. Vassiliev, \textit{Complements of discriminants of
    smooth maps: topology and applications}, Translations of Mathematical
  Monographs, \textbf{98} AMS, Providence, RI, 1992.
\bibitem{vershinin} V. Vershinin, \textit{Braid groups and loop spaces},
  Russian Math. surveys \textbf{ 54} (1999), no. 2, 273--350.
\bibitem{wagner} C. Wagner,
\textit{Symmetric, cyclic, and permutation products of manifolds},
Dissertationes Math. (Rozprawy Mat.)  {\bf 182}  (1980).
\bibitem{wang}  J. Wang, \textit{On the braid groups for $\bbr P^2$ and the
mobius band}, thesis university of rochester, 1997.
\bibitem{welcher} P.J. Welcher, \textit{Symmetric products and the stable
    Hurewicz homomorphism}, Illinois J. Math.  \textbf{24} (1980), no. 4,
  527--544.
\bibitem{zanos} S. Zanos, thesis in progress, university of Lille I.

\bibliographystyle{amsplain}
\bibliography{atest}
\end{thebibliography}
\bibliographystyle{plain[10pt]}

\end{document}